\documentclass[reqno]{amsart}
\usepackage{amsfonts}

\title{ Ricci-quadratic  homogeneous Randers spaces$^*$}
\author{ Shaoqiang Deng \and  Zhiguang  Hu}
\address{School of Mathematical Sciences and LPMC, Nankai University, Tianjin 300071, P R China}
\address{College of Mathematical Sciences, Tianjin Normal University, Tianjin 300387, P.R. China}
\email[Shaoqiang Deng]{dengsq@nankai.edu.cn}
\email[Zhiguang Hu]{nankaitaiji@mail.nankai.edu.cn}

\date{}

\newtheorem{thm}{Theorem}[section]

\newtheorem{lem}[thm]{Lemma}

\theoremstyle{definition}

\thanks{$^*$Project supported by NSFC (no. 10971104, 10671096) and SRFDP of China}
\thanks{Zhiguang Hu is the corresponding author. E-mail address: nankaitaiji@mail.nankai.edu.cn}
\begin{document}

\maketitle

\noindent{\bf Abstract}

\medskip A  Finsler space is called Ricci-quadratic if its Ricci curvature $Ric(x,y)$ is
quadratic in $y$. It is called a Berwald space if its Chern
connection defines a linear connection directly on the underlying
manifold $M$. In this article, we prove that a homogeneous Randers
space is Ricci-quadratic if
 and only if it is of  Berwald type.





\section {\textbf{Introduction}}
Riemann curvature is a central concept in Riemannian geometry which
was introduced by Riemann in 1854. In 1926, Berwald generalized this notion to
Finsler metrics. A Finsler metric is said to be R-quadratic if its
Riemann curvature is quadratic \cite{C-S04}. R-quadratic metrics
were first introduced by Basco and Matsumoto \cite{B-M00}. They form
a rich class of Finsler spaces. For example, all Berwald metrics are
R-quadratic, and some non-Berwald R-quadratic Finsler metrics have
been constructed in \cite{BRS04,L-S09}. There are many interesting works related
to this subject (cf. \cite{S01,NBT07}).

Ricci curvature of a Finsler space is the trace of the Riemann
curvature. A Finsler metric is called Ricci-quadratic if its Ricci curvature $Ric(x,y)$ is quadratic in $y$.
  It is clear that  the notion of Ricci-quadratic metrics
is weaker than that of R-quadratic metrics. It is therefore obvious that any R-quadratic Finsler space
must be Ricci-quadratic, in particular, any Berwald space must be Ricci-quadratic. However,  there are many non-Berwald
spaces which are Ricci-quadratic. In general, it is quite
difficult to characterize Ricci-quadratic metrics. Li and Shen considered
the case of Randers metrics in \cite{L-S09} and obtained a characterization of Ricci-quadratic properties of such spaces, using
some complicated calculations in local coordinate systems. Their results
are rather complicated (see Theorem 3.1 below). In this
paper we consider homogeneous Randers spaces and  prove the following

\noindent \textbf{Main theorem}. {\it  A homogeneous Randers space
is Ricci-quadratic if and only if it is  of  Berwald type}.

We remark here that the range of homogeneous Randers spaces is
rather wide. For example, on any connected Lie group $G$ with Lie
algebra $\mathfrak{g}$, we can identify the tangent space of $G$
at the origin $T_o(G)$ with $\mathfrak{g}$. Given any inner
product $\langle,\rangle$ and a vector $w\in\mathfrak{g}$ with
$\langle w,w\rangle<1$, we can then define a Minkowski norm $F_o$
on $\mathfrak{g}$ by (see \cite{BCS})
$$F_o(u)=\sqrt{\langle u,u\rangle}+\langle w,u\rangle.$$
Then we can extend this Minkowski norm to a left invariant Randers
metric $F$ on $G$ by the left translation of $G$. This method
produces numerous examples of homogeneous Randers spaces. It is
also known that on many coset spaces $G/H$ of a Lie group $G$ with
respect to a non-trivial closed subgroup $H$, there exist
$G$-invariant Randers metrics (see for example \cite{D08}).

\section {\textbf{Preliminaries}}
Let $F$ be a Finsler metric on an $n$-dimensional manifold $M$. We
always assume that $F$ is positive definite, namely, the Hessian
matrix $g_{ij}=g_{ij}(x,y)$ is positive definite, where
   \[ g_{ij}(x,y):= \frac{1}{2}[F^2]_{y^iy^j}(x,y), \quad y\in
      T_xM-\{0\}.  \]
On a standard local coordinate system $(x^1,x^2,\cdots, x^n, y^1, y^2,\cdots, y^n)$, the  geodesics of $F$ are characterized by the
following system of equations:
    \[  \frac{d^2x^i}{dt^2}+2G^i(x,\frac{dx}{dt})=0, \]
where $G^i=G^i(x,y)$ are called the geodesic coefficients of $F$,
which are given by
     \[ G^i=\frac{1}{4}\Big\{ [F^2]_{x^my^l}y^m-[F^2]_{x^l}\Big\}. \]
It is clear that if $F$ is Riemannian, then $G^i$ are quadratic in
$y$. For a general Finsler metric, $G^i$ is very complicated and is not quadratic in $y$.  When $G^i$ are quadratic in $y$, we call $F$  a Berwald
metric. Every Riemannian metric is a Berwald metric, but the converse is not true. In fact, one can construct many examples of non-Riemannian Berwald metrics.
The local structure of Berwald spaces was determined by Z. I. Szab\'{o} in \cite{SZ}.

A Finsler metric of the form $F=\alpha+\beta$, where $\alpha$ is a Riemannian metric and $\beta$ is a $1$-form on $M$ whose length with respect to $\alpha$ is everywhere less than $1$, is called a Randers metric. This kind of metrics was introduced by G. Randers in 1941 (\cite{RA}), in his study of general relativity. A Randers metric  $F=\alpha+\beta$ is a Berwald metric if and only if the form $\beta$ is parallel with respect to $\alpha$. This is an important result in the field of Finsler geometry due to the contributions of many mathematicians, see \cite{BCS} for an account of the history of this result.

Let $y$ be a non-zero vector in $T_x(M)$. The Riemann curvature $R_y = R^i_{\ k} \frac{\partial}{\partial
x^i}\bigotimes dx^k$ is defined by
  \[ R^i_{\ k} := 2[G^i]_{x^k}-[G^i]_{x^my^k}y^m +
       2G^m[G^i]_{y^my^k}-[G^i]_{y^m}[G^m]_{y^k}. \]
It defines a linear transformation on $T_x(M)$. The trace of this linear transformation is denoted by $Ric(x,y)$ and is called the Ricci curvature of $F$.
A Finsler metric is called Ricci quadratic if $Ric(x,y)$ is quadratic in $y$.

\section{The Levi-Civita connection of homogeneous spaces}

In this section we shall use  Killing vector fields to present
some formulas about the  Levi-Civita connection of homogeneous Riemannian manifolds. We follow
the method used by the first author in \cite{D08}.

Let $(G/H,\alpha)$ be a homogeneous Riemannian manifold. Then the
Lie algebra of $G$ has a decomposition $\mathfrak{g}=
\mathfrak{h}+\mathfrak{m}$, where $\mathfrak{h}$ is the Lie
algebra of $H$ and
$Ad(h)(\mathfrak{m})\subset\mathfrak{m},\forall{h}\in{H}$. We
identify $\mathfrak{m}$ with the tangent space $T_{o}(G/H)$ of the
origin $o=H$. We shall use the notation $\langle,\rangle$ to
denote  the Riemannian metric on the manifold as well as its
restriction to $\mathfrak{m}$.  Note that it is an $AdH$-invariant
inner product on $\mathfrak{m}$. Hence we have
$$\langle [x,u],v\rangle+\langle[x,v],u\rangle=0,\ \ \forall
x\in\mathfrak{h}, \forall u,v\in\mathfrak{m},$$ which is equivalent
to
$$\langle[x,u],u\rangle=0, \ \ \forall x\in\mathfrak{h},
\forall u\in\mathfrak{m}.$$

Given $v\in \mathfrak{g}$, we can define the fundamental vector field
$\hat{v}$ generated by $v$, i.e.,
$$\hat{v}_{gH}=\frac{d}{dt}\exp (tv)gH|_{t=0},\indent \forall g\in G.$$
Since the  one-parameter transformation group $\exp tv$ on $G/H$
consists of isometries, $\hat{v}$ is a Killing vector field.

Let $\hat{X},\hat{Y},\hat{Z}$ be Killing vector fields on $G/H$ and
$U,V$ be arbitrary smooth vector fields on $G/H$. Then we have (\cite{B87}, page 40,182,183)
\begin{align}
&[\hat{X},\hat{Y}]=-[X,Y]^{\hat{}},\label{kil0}\\
&\hat{X}\langle
U,V\rangle=\langle[\hat{X},U],V\rangle+\langle[\hat{X},V],U\rangle,\label{kil1}\\
&\langle\nabla_{\hat{X}}\hat{Y},\hat{Z}\rangle
=-\frac{1}{2}\left(\langle[X,Y]^{\hat{}},\hat{Z}\rangle+\langle[X,Z]^{\hat{}},\hat{Y}\rangle
+\langle[Y,Z]^{\hat{}},\hat{X}\rangle\right).\label{kil2}
\end{align}
We only need to prove (2.2). In fact, by (c) of Theorem 1.81 of \cite{B87}, we have:
\begin{align*}
\hat{X}\langle
U,V\rangle &=\langle\nabla_{\hat{X}}U, V\rangle+\langle U, \nabla_{\hat{X}}V\rangle\\
&=\langle \nabla_U \hat{X}, V\rangle+[\hat{X}, U], V\rangle+\langle U, \nabla_V\hat{X}\rangle+\langle U, [\hat{X}, V]\rangle\\
&=\langle[\hat{X},U],V\rangle+\langle[\hat{X},V],U\rangle.
\end{align*}

Let $u_1,u_2,\cdots,u_n$ be an orthonormal basis of $\mathfrak{m}$
with respect to $\langle,\rangle$. We extend it to  a basis
$u_1,u_2,\cdots,u_m$ of $\mathfrak{g}$. By \cite{H78}, there exists
a local coordinate system on a neighborhood $V$ of $o$, which is
defined by the mapping
$$(\exp(x^{1}u_1)\exp (x^2u_2)\cdots \exp (x^nu_n))H\rightarrow (x^1,x^2,\cdots,x^n).$$

Let $gH=(x^1,x^2,\cdots,x^n)\in U$. Then
$$\left.\frac{\partial}{\partial{x^i}}\right|_{gH}=f_i^{a}\hat{u_a}|_{gH},$$
where
 \begin{equation} f^{a}_iu_a=e^{x^1adu_1}\cdots
e^{x^{i-1}adu_{i-1}}(u_i).\label{fff}
\end{equation}

\textbf{\emph{Remark}}  In the following, the indices
$a,b,c,\cdots$ range from $1$ to $m$, the indices $i,j,k,\cdots$
range from $1$ to $n$ and the indices $\lambda, \mu,\cdots$ range
from $n+1$ to $m$.

Let $\Gamma_{ij}^l$ be the Christoffel symbols in the coordinate
system, i.e., $$\nabla_{\frac{\partial}{\partial
x^i}}\frac{\partial}{\partial
x^j}=\Gamma_{ij}^k\frac{\partial}{\partial x^k}.$$ Then
\begin{equation}
\Gamma_{ij}^l\frac{\partial}{\partial
x^l}=\nabla_{\frac{\partial}{\partial x^i}}\frac{\partial}{\partial
x^j}=\frac{\partial f_j^{a}}{\partial
x^i}\hat{u_a}+f_i^{b}f_j^{a}\nabla_{\hat{u_b}}\hat{u_a}.\label{gam}
\end{equation}

\noindent From \eqref{fff}, we see that $f_i^{a}$ are functions of
$x^1,\cdots,x^{i-1}$. Thus $$\frac{\partial f_j^{a}}{\partial
x^i}=0,\ \ i\geq j.$$ Therefore \eqref{gam} gives
\begin{equation*}
\Gamma_{ij}^l\frac{\partial}{\partial
x^l}=f_i^{b}f_j^{a}\nabla_{\hat{u_b}}\hat{u_a},\ \ i\geq j.
\end{equation*}

\noindent Differentiating the above equation with respect to $x_k$,
we get
\begin{align}
\frac{\partial \Gamma_{ij}^l}{\partial
x^k}\frac{\partial}{\partial
x^l}+\Gamma_{ij}^s\Gamma_{ks}^l\frac{\partial}{\partial
x^l}=\frac{\partial (f_i^{b}f_j^{a})}{\partial
x^k}\nabla_{\hat{u_b}}\hat{u_a}+f_i^{b}f_j^{a}f_k^{c}\nabla_{\hat{u_c}}\nabla_{\hat{u_b}}\hat{u_a},\
\ i\geq j.\label{pgam}
\end{align}

\noindent Differentiating \eqref{fff} with respect to $x_k$ and
letting $(x^1,\cdots,x^n)\rightarrow0$, we obtain
\begin{equation*}
\frac{\partial f_i^{a}}{\partial x^k}(0)=f(k,i)C_{ki}^{a},
\end{equation*}
where $C_{ab}^c$ are the structure constants of $\mathfrak{g}$
which are defined by $[u_a, u_b]=C_{ab}^cu_c$ and $f(k,l)$ are
defined by
$$f(k,i):=\mathbf{}
     \begin{cases} 1,\ &k<i,\\ 0&k\geq i.
     \end{cases}$$

Considering the value at the origin $o$ we get
\begin{lem}
\begin{align}
&\Gamma_{ij}^l(o)=f(i,j)C_{ij}^l+\langle\nabla_{\hat{u_i}}\hat{u_j},\hat{u_l}\rangle\label{gam1},\\
&\begin{aligned} \left.\frac{\partial \Gamma_{ij}^l}{\partial
x^k}\right|_o=&-\Gamma_{ij}^s(\Gamma_{ks}^l+\langle\nabla_{\hat{u_k}}\hat{u_l},\hat{u_s}\rangle)
+f(k,j)C_{kj}^a\langle\nabla_{\hat{u_i}}\hat{u_a},\hat{u_l}\rangle\\
&+f(k,i)C_{ki}^s\langle\nabla_{\hat{u_s}}\hat{u_j},\hat{u_l}\rangle
+\hat{u_k}\langle\nabla_{\hat{u_i}}\hat{u_j},\hat{u_l}\rangle,\ \
i\geq j\label{gam2}.
\end{aligned}
\end{align}
\end{lem}

\textbf{\emph{Proof}}. From \eqref{fff} we know that  $f_i^a(0)=\delta
_i^a$ and that $\frac{\partial}{\partial x^k}|_o=\hat{u_k}|_o=u_k$. Thus,
by \eqref{gam}
$$\Gamma_{ij}^l(o)=\langle\Gamma_{ij}^s\frac{\partial}{\partial
x^s}, \hat{u_l}\rangle|_o=\langle\frac{\partial f_j^{a}}{\partial
x^i}\hat{u_a}+f_i^{b}f_j^{a}\nabla_{\hat{u_b}}\hat{u_a},\hat{u_l}\rangle|_o.$$
\eqref{gam1} is obtained from the above equation. By \eqref{pgam} we
get
\begin{align*}
\left.\frac{\partial \Gamma_{ij}^l}{\partial
x^k}\right|_o=&-\Gamma_{ij}^s\Gamma_{ks}^l
+f(k,j)C_{kj}^a\langle\nabla_{\hat{u_i}}\hat{u_a},\hat{u_l}\rangle\\
&+f(k,i)C_{ki}^a\langle\nabla_{\hat{u_a}}\hat{u_j},\hat{u_l}\rangle
+\langle\nabla_{\hat{u_k}}\nabla_{\hat{u_i}}\hat{u_j},\hat{u_l}\rangle,\
\ i\geq j.
\end{align*}
And we know that at the origin
$$\langle\nabla_{\hat{u_k}}\nabla_{\hat{u_i}}\hat{u_j},\hat{u_l}\rangle=
\hat{u_k}\langle\nabla_{\hat{u_i}}\hat{u_j},\hat{u_l}\rangle
-\langle\nabla_{\hat{u_i}}\hat{u_j},\nabla_{\hat{u_k}}\hat{u_l}\rangle,$$
$$f(k,i)C_{ki}^a\langle\nabla_{\hat{u_a}}\hat{u_j},\hat{u_l}\rangle=
f(k,i)C_{ki}^s\langle\nabla_{\hat{u_s}}\hat{u_j},\hat{u_l}\rangle,$$
and
$$\langle\nabla_{\hat{u_i}}\hat{u_j},\nabla_{\hat{u_k}}\hat{u_l}\rangle
=\Gamma_{ij}^s\langle\nabla_{\hat{u_k}}\hat{u_l},\hat{u_s}\rangle, \
\ i\geq j.$$ From the above four equations \eqref{gam2} is obtained.
\ \ $\Box$

We will also need the following
\begin{lem}
For $u_i,u_j,u_k,u_l\in\mathfrak{m},u_\lambda\in\mathfrak{h}$, we
have
\begin{align}
&\langle\nabla_{\hat{u_i}}\hat{u_j},\hat{u_l}\rangle|_o=-\frac{1}{2}(C_{ij}^l+C_{il}^j+C_{jl}^i)\label{lem20},\\
&\langle\nabla_{\hat{u_i}}\hat{u_\lambda},\hat{u_j}\rangle|_o=\langle[u_j,
u_\lambda]_{\mathfrak{m}},u_i\rangle=C_{j\lambda}^i\label{lem21},\\
&\begin{aligned}
\hat{u_k}\langle\nabla_{\hat{u_i}}\hat{u_j},\hat{u_l}\rangle|_o=\frac{1}{2}\left(C_{ka}^lC_{ij}^a+C_{ka}^jC_{il}^a
+C_{ka}^iC_{jl}^a+C_{ij}^sC_{kl}^t\delta_{st}+C_{il}^sC_{kj}^t\delta_{st}+C_{jl}^sC_{ki}^t\delta_{st}\right),
\end{aligned}\label{lem22}
\end{align}
where $[v_i,v_j]_{\mathfrak{m}}$ denotes the projection of
$[v_i,v_j]$ to $\mathfrak{m}$.
\end{lem}
\textbf{\emph{Proof}}. First, \eqref{lem20} is an alternative
formulation of \eqref{kil2} at the origin in terms of the structure
constants. Using the invariance of $adu_\lambda$, \eqref{lem21} can
also be deduced from
 \eqref{kil2}. Finally,  by \eqref{kil2} we have
 $$\hat{u_k}\langle\nabla_{\hat{u_i}}\hat{u_j},\hat{u_l}\rangle=-\frac{1}{2}\hat{u_k}
\left(\langle[u_i,u_j]^{\hat{}},\hat{u_l}\rangle+\langle[u_i,u_l]^{\hat{}},\hat{u_j}\rangle
+\langle[u_j,u_l]^{\hat{}},\hat{u_i}\rangle\right).$$ Considering the value at
the origin $o$ and taking into account  \eqref{kil1} and \eqref{kil0},  we can deduce from the above equation that
\begin{align*}
\hat{u_k}\langle\nabla_{\hat{u_i}}\hat{u_j},\hat{u_l}\rangle|_o=&\frac{1}{2}\Big(
\langle\big[u_k,[u_i,u_j]\big]_{\mathfrak{m}},u_l\rangle+\langle\big[u_k,[u_i,u_l]\big]_{\mathfrak{m}},u_j\rangle\\
&+\langle\big[u_k,[u_j,u_l]\big]_{\mathfrak{m}},u_i\rangle
+\langle[u_i,u_j]_{\mathfrak{m}},[u_k,u_l]_{\mathfrak{m}}\rangle\\
&+\langle[u_i,u_l]_{\mathfrak{m}},[u_k,u_j]_{\mathfrak{m}}\rangle
+\langle[u_j,u_l]_{\mathfrak{m}},[u_k,u_i]_{\mathfrak{m}}\rangle\Big),
\end{align*}
from which we  get \eqref{lem22}. \ \ $\Box$

By the above two lemmas,  at the origin $o$ we have
\begin{align*}
\Gamma_{ni}^j-\Gamma_{nj}^i=\langle\nabla_{\hat{u_n}}\hat{u_i},\hat{u_j}\rangle
-\langle\nabla_{\hat{u_n}}\hat{u_j},\hat{u_i}\rangle =&C_{ji}^n.
\end{align*}

\section{Ricci-quadratic homogeneous
Randers spaces} In this section we will recall some basic notations
about Randers spaces. Let
$$F=\alpha+\beta=\sqrt{a_{ij}(x)y^iy^j}+b_i(x)y^i$$ be a Randers
metric. Let $\nabla\beta=b_{i|j}y^idx^j$ denote the covariant
derivative of $\beta$ with respect to $\alpha$.
\begin{align*}
r_{ij}:=\frac{1}{2}(b_{i|j}+b_{j|i}),\qquad
s_{ij}:=\frac{1}{2}(b_{i|j}-b_{j|i}), \qquad s_j:=b^{i}s_{ij},
\qquad t_j:=s_m{s^m}_j.
\end{align*}
We use $a_{ij}$ to raise and lower the indices of tensors defined by
$b_i$ and $b_{i|j}$. The index ``0'' means the contraction with
$y^i$. For example, $s_0=s_iy^i$ and $r_{00}=r_{ij}y^iy^j$, etc. For
Ricci-quadratic metrics on Randers spaces, we have the following

\begin{thm}\cite{L-S09}
Let $F=\alpha+\beta$ be a Randers metric on an $n$-dimensional
manifold. Then it is  Ricci-quadratic if and only if
\begin{align}
\label{ric_1}r_{00}+2s_0\beta=2\tilde{c}(\alpha^2-\beta^2),\\
\label{ric_2}{s^k}_{0|k}=(n-1)A_0,
\end{align}
where $\tilde{c}=\tilde{c}(x)$ is a scalar function and
$A_{k}:=2\tilde{c}s_{k}+\tilde{c}^2b_k+t_k+\frac{1}{2}\tilde{c}_k$, here $\tilde{c}_k=\frac{\partial\tilde{c}}{\partial x^k}$.
\end{thm}

Now we consider homogeneous Randers spaces. Let  $(G/H, \alpha)$
and ${\frak m}$ be as above.  If $W$ is a $G$-invariant vector
field on $G/H$, then the restriction of $W$ to $T_o(G/H)$ must be
fixed by the isotropy action of $H$. Under the  identification of
$T_o(G/H)$ with ${\frak m}$, $W$ corresponds to a vector $w\in
{\frak m}$ which is fixed by $Ad(H)$. On the other hand, if $w\in
{\frak m}$ is fixed by $Ad(H)$, then we can define a vector field
$W$ on $G/H$ by $W|_{gH}=\frac{d}{dt}(g\exp (tw)H)|_{t=0}$.
Therefore, $G$-invariant vector fields on $G/H$ are one-to-one
corresponding to vectors in ${\frak m}$ fixed by $Ad(H)$. Note
that a Randers space $F=\alpha+\beta$ is $G$-invariant if and only
if $\alpha$ and $\beta$ are both invariant under $G$. Through
$\alpha$,  $\beta$ corresponds to a vector field $\widetilde{U}$
which is  invariant under $G$ and satisfying $\alpha
(\widetilde{U})<1$ everywhere.
 This implies that  there is a one-to-one correspondence
between the invariant Randers metrics on $G/H$ with the underlying
Riemannian metric and the set $$V=\{u\in\mathfrak{m}|Ad(h)u=u,
\langle u,u\rangle<1,\quad\forall h\in H\}.$$ Also note that in
this case the length $c$ of $\beta$ (or $\widetilde{U}$) is
constant.

Let $(G/H,F)$ be a homogeneous Randers space and
$(U,(x^1,\cdots,x^n))$ be a local coordinate system as in Section 2.
We suppose the vector field $\widetilde{U}$ which corresponds to the
invariant 1-form $\beta$ corresponds to $u=cu_n(c<1)$ under the
Riemannian metric $\alpha$. Thus
\begin{align*}
\widetilde{U}|_{gH}=&\frac{d}{dt}g\exp(tu)H|_{t=0}\\
=&\frac{d}{dt}(\exp x^1u_1\exp x^2u_2\cdots \exp(x^n+ct)u_n)H|_{t=0}\\
=&c\frac{\partial}{\partial x^n}|_{gH}.
\end{align*}
Then we have the following (see\cite{D08}):
\begin{align}\begin{split}
&b_i=\beta(\frac{\partial}{\partial
x^i})=\langle\widetilde{U},\frac{\partial}{\partial
x^i}\rangle=c\langle\frac{\partial}{\partial
x^n},\frac{\partial}{\partial
x^i}\rangle=ca_{ni},\\
&\frac{\partial b_i}{\partial x^j}=c\frac{\partial a_{ni}}{\partial
x^j}=c(\Gamma_{nj}^{k}a_{ki}+\Gamma_{ji}^{k}a_{kn}),\\
&b_{i|j}=\frac{\partial b_i}{\partial
x^j}-b_{l}\Gamma_{ij}^{l}=c\Gamma_{nj}^{k}a_{ki},\\
&r_{ij}=\frac{1}{2}(b_{i|j}+b_{j|i})=\frac{c}{2}(\Gamma_{nj}^{k}a_{ki}+\Gamma_{ni}^{k}a_{kj}),\\
&s_{ij}=\frac{1}{2}(b_{i|j}-b_{j|i})=\frac{c}{2}(\Gamma_{nj}^{k}a_{ki}-\Gamma_{ni}^{k}a_{kj}),\\
&s_{j}=b^{i}s_{ij}=a^{il}b_{l}s_{ij}=cs_{nj}.
\end{split}\label{bij}\end{align}

\begin{lem}
Let $(G/H,F)$ be a homogeneous Randers space and $\beta$ correspond
to $u$. Then \eqref{ric_1} implies that
\begin{equation} \langle[y,u]_{\mathfrak{m}},y\rangle=0,\quad
\forall y\in\mathfrak{m}. \label{qric_1}\end{equation}
\end{lem}
\textbf{\emph{Proof}}. Considering the value at $o$, by
\eqref{bij},  \eqref{gam1} and \eqref{lem20}, we have
\begin{align*}
(r_{00}+2s_0\beta)|_{o}&=c\Gamma_{n0}^0+2c\frac{c}{2}(\Gamma_{n0}^{n}-\Gamma_{nn}^{0})\langle
u,y\rangle\\
&=cC_{0n}^0+c^2C_{n0}^n\langle
u,y\rangle\\
&=\langle[y,u]_{\mathfrak{m}},y-\langle u,y\rangle u\rangle.
\end{align*}
It is obvious that
\begin{align*}
2\tilde{c}(\alpha^2-\beta^2)|_{o}=2\tilde{c}(o)(\langle
y,y\rangle-\langle u,y\rangle^2).
\end{align*}
Plugging the above two equations into \eqref{ric_1}, we get that at
$o$
\begin{align*}
\Big\langle[y,u]_{\mathfrak{m}},y-\langle u,y\rangle
u\Big\rangle=2\tilde{c}(o)\left(\langle y,y\rangle-\langle
u,y\rangle^2\right).
\end{align*}
Setting $y=u$ and taking into account the fact that $\langle
u,u\rangle<1$ , we get
$$\tilde{c}(o)=0.$$
Thus $$\Big\langle [y,u]_{\mathfrak{m}}, y-\langle u,y\rangle
u\Big\rangle=0.$$ Replacing $y$ by $y+u$ in the above equation
yields $$\Big\langle [y,u]_{\mathfrak{m}}, u+y-\langle u,u+y\rangle
u\Big\rangle=0.$$
From the above two equations and the fact that $\langle
u,u\rangle<1$ we deduce that
$$\langle[y,u]_{\mathfrak{m}},u\rangle=0,\quad
\langle[y,u]_{\mathfrak{m}},y\rangle=0.$$ This proves the lemma. \ \
$\Box$

Note that in the above  lemma we  also have
\begin{align}
&C_{ni}^j+C_{nj}^i=0=C_{ni}^n,\label{anti}\\
&s_{i}(o)=cs_{ni}=\frac{c^2}{2}(\Gamma_{ni}^{n}-\Gamma_{nn}^{i})
=\frac{1}{2}\langle[u,u_{i}]_{\mathfrak{m}},u\rangle=0,\nonumber\\
&t_i(o)=s_m{s^m}_i=0.\nonumber
\end{align}

\section{Proof of the main theorem}
Before the proof, we still need to perform some complicated computation.

First we have

\begin{align}
\begin{split}
\frac{\partial b_{i|j}}{\partial x^k}&=\frac{\partial
\Gamma_{nj}^l}{\partial x^k}
+\Gamma_{nj}^s(\Gamma_{ks}^l+\Gamma_{kl}^t\delta_{ts})\\
&=f(k,l)C_{kl}^s\Gamma_{nj}^t\delta_{ts}+f(k,j)C_{kj}^a\langle\nabla_{\hat{u_n}}\hat{u_a},\hat{u_l}\rangle\\
&\ \ +C_{kn}^s\langle\nabla_{\hat{u_s}}\hat{u_j},\hat{u_l}\rangle
+\hat{u_k}\langle\nabla_{\hat{u_n}}\hat{u_j},\hat{u_l}\rangle.
\end{split}\label{ups}
\end{align}


\begin{lem}
Let $(G/H,F)$ be a homogeneous Randers space and $\beta$
correspond to $u$ which satisfies \eqref{qric_1}. Then
\begin{equation} \frac{\partial b_{0|0}}{\partial x^0}=0. \label{000}\end{equation}
\end{lem}

\textbf{\emph{Proof}}. By the Lemma 3.2 
we have the following computations
\begin{align}
\begin{split}
\frac{\partial b_{0|0}}{\partial
x^0}&=f(i,j)C_{ij}^sy^iy^j\Gamma_{n0}^t\delta_{ts}+f(i,j)C_{ij}^sy^iy^j\langle\nabla_{\hat{u_n}}\hat{u_s},\hat{u_0}\rangle
+f(i,j)C_{ij}^\lambda y^iy^j\langle\nabla_{\hat{u_n}}\hat{u_\lambda},\hat{u_0}\rangle\\
&\ \ +C_{0n}^s\langle\nabla_{\hat{u_s}}\hat{u_0},\hat{u_0}\rangle
+\hat{u_0}\langle\nabla_{\hat{u_n}}\hat{u_0},\hat{u_0}\rangle\\
&=f(i,j)C_{ij}^sy^iy^j(\langle\nabla_{\hat{u_n}}\hat{u_0},\hat{u_s}\rangle+\langle\nabla_{\hat{u_n}}\hat{u_s},\hat{u_0}\rangle)
+f(0,0)C_{00}^\lambda y^iy^jC_{\lambda n}^0-C_{0n}^sC_{s0}^0
+C_{0a}^0C_{n0}^a\\
&=f(i,j)C_{ij}^sy^iy^j(C_{n0}^s+C_{ns}^0)+f(i,j)C_{ij}^\lambda
y^iy^j C_{\lambda
n}^0+C_{0\lambda}^0C_{n0}^\lambda\\
&=f(i,j)C_{ij}^\lambda y^iy^jC_{\lambda n}^0\\
&=0.
\end{split}
\end{align}
In above we have used the fact $C_{n\lambda}^i=0$.\ \ $\Box$

 Further,
we also have the following
\begin{lem}
Let $(G/H,F)$ be a homogeneous Randers space and $\beta$
correspond to $u$ which satisfies \eqref{qric_1}. Then
\begin{equation}\frac{\partial{s_{k0}}}{\partial{x^i}}|_o=
\frac{c}{2}\big(f(i,k)C_{ik}^sC_{s0}^n+f(i,0)C_{i0}^sC_{ks}^n+C_{0k}^sC_{is}^n\big).\label{sijd}
\end{equation}
\end{lem}

\textbf{\emph{Proof}}. By \eqref{ups}, we get
\begin{align*}
\frac{2}{c}\frac{\partial{s_{k0}}}{\partial{x^i}}|_o
=&\frac{\partial b_{k|0}}{\partial x^i}-\frac{\partial b_{0|k}}{\partial x^i}\\
=&f(i,k)C_{ik}^s\Gamma_{n0}^t\delta_{st}+f(i,0)C_{i0}^a\langle\nabla_{\hat{u_n}}\hat{u_a},\hat{u_k}\rangle+
C_{in}^s\langle\nabla_{\hat{u_s}}\hat{u_0},\hat{u_k}\rangle\\&\ \
-(f(i,0)C_{i0}^s\Gamma_{nk}^t\delta_{st}+
f(i,k)C_{ik}^a\langle\nabla_{\hat{u_n}}\hat{u_a},\hat{u_0}\rangle
+C_{in}^s\langle\nabla_{\hat{u_s}}\hat{u_k},\hat{u_0}\rangle)\\
&\ \ \ +\hat{u_i}\langle\nabla_{\hat{u_n}}\hat{u_0},\hat{u_k}\rangle
-\hat{u_i}\langle\nabla_{\hat{u_n}}\hat{u_k},\hat{u_0}\rangle\\
=&f(i,0)(C_{i0}^a\langle\nabla_{\hat{u_n}}\hat{u_a},\hat{u_k}\rangle
-C_{i0}^s\Gamma_{nk}^t\delta_{st})-f(i,k)(C_{ik}^a\langle\nabla_{\hat{u_n}}\hat{u_a},\hat{u_0}\rangle
-C_{ik}^s\Gamma_{n0}^t\delta_{st})\\
&\ \
+C_{in}^s(\langle\nabla_{\hat{u_s}}\hat{u_0},\hat{u_k}\rangle-\langle\nabla_{\hat{u_s}}\hat{u_k},\hat{u_0}\rangle)
+\hat{u_i}\langle\nabla_{\hat{u_n}}\hat{u_0},\hat{u_k}\rangle
-\hat{u_i}\langle\nabla_{\hat{u_n}}\hat{u_k},\hat{u_0}\rangle
\end{align*}
By \eqref{gam1}, \eqref{lem20}, we have
\begin{align*}
\langle\nabla_{\hat{u_s}}\hat{u_0},\hat{u_k}\rangle-\langle\nabla_{\hat{u_s}}\hat{u_k},\hat{u_0}\rangle
=C_{k0}^s.
\end{align*}
Then by \eqref{lem21}, we get
\begin{align*}
C_{i0}^a\langle\nabla_{\hat{u_n}}\hat{u_a},\hat{u_k}\rangle
-C_{i0}^s\Gamma_{nk}^t\delta_{st}&=C_{i0}^s
(\langle\nabla_{\hat{u_n}}\hat{u_s},\hat{u_k}\rangle-\langle\nabla_{\hat{u_n}}\hat{u_k},\hat{u_s}\rangle)
+C_{i0}^{\lambda}\langle\nabla_{\hat{u_n}}\hat{u_{\lambda}},\hat{u_k}\rangle\\
&=C_{i0}^sC_{ks}^n.
\end{align*}
From \eqref{lem22} we easily get
\begin{align*}
\hat{u_i}\langle\nabla_{\hat{u_n}}\hat{u_0},\hat{u_k}\rangle
-\hat{u_i}\langle\nabla_{\hat{u_n}}\hat{u_k},\hat{u_0}\rangle=C_{is}^nC_{0k}^s+C_{0k}^sC_{in}^t\delta_{st}.
\end{align*}
Combining the above four equations, we obtain \eqref{sijd}. \ \
$\Box$

\textbf{\emph{Proof of the main theorem}}. Let $G/H$ be a
homogeneous Randers space. Taking the local coordinate system as in
Section 3, we have seen that \eqref{qric_1} holds. In particular, we
have $C_{ni}^n=0$. 
By \eqref{sijd} we have
\begin{align*}
\left.\frac{\partial{s_{n0}}}{\partial{x^0}}\right|_o=\frac{c}{2}\Big(f(0,n)C_{0n}^sC_{s0}^n
+f(0,0)C_{00}^sC_{ns}^n+C_{0n}^sC_{0s}^n\Big)=0.
\end{align*}
Differentiating \eqref{ric_1} and taking into account the fact that
$\tilde{c}(o)=0=s_0(o)$,  we deduce from \eqref{bij} that
\begin{align*}
2\tilde{c_{0}}(\alpha^2-\beta^2)|_{o}=c\frac{\partial
b_{0|0}}{\partial x^0}+2c\beta
\frac{\partial{s_{n0}}}{\partial{x^{0}}}=0.
\end{align*}
Thus
$$\tilde{c_0}=0,$$
which yields
$$A_0(o)=\frac{1}{2}\tilde{c}_0(o)=0.$$

On the other hand, we have
\begin{align*}
s_{ij}(o)=\frac{c}{2}(\Gamma_{nj}^i-\Gamma_{ni}^j)=\frac{c}{2}C_{ij}^n.
\end{align*}
Thus{\small
\begin{align*}
{s^{k}}_{0|k}(o)&=\sum_ks_{k0|k}=\sum_{k,l}(\frac{\partial{s_{k0}}}{\partial{x^k}}-\Gamma_{kk}^ls_{l0}-\Gamma_{0k}^ls_{kl})\\
&=\frac{c}{2}\sum_{k=1}^n\left(f(k,0)-1\right)C_{k0}^sC_{ks}^n
+\sum_{1\leq k,l\leq
n}\frac{c}{2}\left(C_{l0}^nC_{kl}^k-C_{kl}^n\left(f(0,k)C_{0k}^l-\frac{1}{2}(C_{0k}^l+C_{0l}^{k}+C_{kl}^0)\right)\right)\\
&= \sum_{1\leq k,l\leq
n}\frac{c}{2}\left(C_{l0}^nC_{kl}^k+\frac{1}{2}C_{kl}^n(C_{0k}^l+C_{0l}^{k}+C_{kl}^0)\right).
\end{align*}}
Therefore \eqref{ric_2} reduces to{\small
$$\sum_{1\leq k,l\leq
n}\frac{c}{2}\left(C_{l0}^nC_{kl}^k+\frac{1}{2}C_{kl}^n(C_{0k}^l+C_{0l}^{k}+C_{kl}^0)\right)=0.$$}
Set $y=u$ in the above equation. Then by \eqref{qric_1} we have
$$C_{kl}^n=0,\qquad k,l=1,\cdots,n,$$
 i.e.,
\begin{align}\label{qric_2}\langle[u_k,u_l]_{\mathfrak{m}},u\rangle=0,\qquad
k,l=1,\cdots,n. \end{align}
It is proved in \cite{D08} that a homogeneous Randers space is of the Berwald type if and only if  \eqref{qric_1} and
\eqref{qric_2} hold. This  completes the proof.
\ \ $\Box$

\end{document}